\documentclass[11pt]{article}
\usepackage{mathrsfs}
\usepackage{amsfonts}
\usepackage{mathrsfs}
\usepackage{amsmath}
\usepackage{amssymb}
\usepackage{graphicx}
\graphicspath{{figure/}}
\usepackage{epic}
\usepackage{amsthm,latexsym,CJK}
\renewcommand{\paragraph}{\roman{paragraph}}
\def \n{\nabla}
\newtheorem{theorem}{Theorem}[section]
\newtheorem{remark}[theorem]{Remark}
\newtheorem{lemma}[theorem]{Lemma}
\newtheorem{corollary}[theorem]{Corollary}

\begin{document}
\title{\bf Constrained matrix Li-Yau-Hamilton estimates \\ on K\"{a}hler manifolds\thanks{
 Email: renx@cumt.edu.cn, ys@cumt.edu.cn, ljshen@cumt.edu.cn, gyzhang@cumt.edu.cn}}
\author{Xin-An Ren\ \  Sha Yao \ \   Li-Ju Shen \ \  Guang-Ying Zhang\\
\small Department of Mathematics, China University of Mining
and Technology, Xuzhou 221116, China.}

\date{}\maketitle\renewcommand{\theequation}{\thesection.\arabic{equation}}

\begin{abstract}
We derive an interpolation version of constrained matrix Li-Yau-Hamilton estimate on K\"{a}hler manifolds. As a result, we first get a constrained matrix Li-Yau-Hamilton estimate for heat equation on a K\"{a}hler manifold with fixed K\"{a}hler metric. Secondly, we get a corresponding estimate for forward conjugate heat equation on K\"{a}hler manifolds with time dependent K\"{a}hler metrics evolving by K\"{a}hler-Ricci flow.
\end{abstract}
\section{Introduction}
In the seminar paper \cite{LY}, Li and Yau first proved a gradient estimate for the positive solution of heat equation
 \begin{equation}
\frac{\partial}{\partial t}u=\Delta u
 \end{equation}
 on a Riemannian manifold and derived the classical Harnack inequality by integrating the gradient estimate along the paths in space-time. This estimate has tremendous impact in the past twenty years. On the one hand, similar technique was employed by Hamilton in study of  Harnack inequalities for Ricci flow in \cite{H3,H5} and mean curvature flow in \cite{H7}. So this class of Harnack inequalities was called Li-Yau-Hamilton estimate in \cite{NT1}. Similar Li-Yau-Hamilton estimates for Gauss curvature flow and Yamabe flow were proved by Chow in \cite{Ch1,Ch3}. Meanwhile, the estimate for K\"{a}hler-Ricci flow was proved by Cao in \cite{C2}. Later in \cite{A}, the corresponding result was established for a general class of hypersurface flows by Andrews. Moreover, a geometric interpretation of Hamilton quantities in the Li-Yau-Hamilton estimate for Ricci flow was given by Chow and Chu in \cite{CC1}. Using this geometric approach, Chow and Knopf gave a new Li-Yau-Hamilton estimate for Ricci flow in \cite{CK}.

  On the other hand, Hamilton extended Li-Yau's gradient estimate of heat equation to the full matrix version in \cite{H4}. Precisely, he proved that a positive solution of heat equation
on a closed Riemannian manifold with nonnegative sectional curvature and parallel Ricci curvature satisfies the estimate
\begin{equation}
\nabla_i\nabla_j\ln u+\frac{1}{2t}g_{ij}\geq 0.
\end{equation}
Moreover, he derived several important monotonicity formulae out of the estimate in \cite{H6}.
 The assumption that $(M,g)$ has parallel Ricci curvature in Hamilton's estimate is rather restrictive and can be removed if the manifold is K\"{a}hler. This was observed by Cao and Ni for the K\"{a}hler manifold with fixed metric in \cite{CN} and by Chow and Ni for the K\"{a}hler manifold with time dependent metrics evolving by K\"{a}hler-Ricci flow in \cite{N}. These two estimates can be unified by the interpolation consideration originally proposed by Chow in \cite{Ch4}.
\vspace{2mm}\newline
{\bf Theorem A}(Chow-Ni) {\it Assume that $(M,g_{i\bar{j}}(t))$ has nonnegative holomorphic bisectional curvature satisfying $\epsilon$-K\"{a}hler-Ricci flow
 \begin{equation}
 \frac{\partial}{\partial t}g_{i\bar{j}}=-\epsilon R_{i\bar{j}}.
 \end{equation}
 In the case that $M$ is complete noncompact, additionally assume that the bisectional curvature is bounded. If $u$ is a positive solution of
\begin{equation}
\frac{\partial}{\partial t}u=\Delta u+\epsilon Ru,
\end{equation}
where $R$ is the scalar curvature, then
\begin{equation}
\nabla_i\nabla_{\bar{j}}\ln u+\epsilon R_{i\bar{j}}+\frac{1}{t}g_{i\bar{j}}\geq 0.
\end{equation}}
It should be pointed out that there are also other important works along this direction, see \cite{CC2,N2,N4,NT1,NT2}.

To better understand how the Li-Yau-Hamilton estimate for the Ricci flow can be perturbed or extended, which is important in studying singularities of the Ricci flow, Chow and Hamilton extended the matrix Li-Yau-Hamilton estimate on Riemannian manifolds to the constrained case in \cite{CH}.
\vspace{2mm}\newline
{\bf Theorem B}(Chow-Hamilton) {\it Let $(M,g)$ be a compact Riemannian manifold with nonnegative sectional curvature and parallel Ricci curvature. If $u$ and $v$ are positive solutions of heat equation
\begin{equation}
\frac{\partial}{\partial t}u=\Delta u,\hspace{4mm}\frac{\partial}{\partial t}v=\Delta v,
\end{equation}
with $|v|<u$, then
\begin{equation}
\nabla_i\nabla_j\ln u+\frac{1}{2t}g_{ij}> \frac{\n_ih\n_{j}h}{1-h^2},
\end{equation}
where $h=v/u$.}\vspace{2mm}

In this paper, we will extend the Li-Yau-Hamilton estimates due to Cao-Ni and Chow-Ni on K\"{a}hler manifolds to the constrained case.
\begin{theorem}
Let $(M,g_{i\bar{j}}(t))$ be a compact solution of $\epsilon$-K\"{a}hler-Ricci flow
\begin{equation}
 \frac{\partial}{\partial t}g_{i\bar{j}}=-\epsilon R_{i\bar{j}}.
 \end{equation}
 If $u$ and $v$ are solutions of the equation
\begin{equation}
\frac{\partial}{\partial t}u=\Delta u+\epsilon Ru, \hspace{6mm}\frac{\partial}{\partial t}v=\Delta v+\epsilon Rv \label{heatequ}
\end{equation}
with $|v|<u$, then we have
\begin{equation}
\n_i\n_{\bar{j}}\ln u+\frac{1}{t}g_{i\bar{j}}+\epsilon R_{i\bar{j}}> \frac{\n_ih\n_{\bar{j}}h}{1-h^2},
\end{equation}
where $h=v/u$.
\end{theorem}

Set $\epsilon=0$, and we get the constrained Li-Yau-Hamilton estimate on K\"{a}hler manifolds.
\begin{corollary}
Let $M$ be a compact K\"{a}hler manifold with nonnegative holomorphic bisectional curvature.
 If $u$ and $v$ are the solutions of the heat equation
\begin{equation}
\frac{\partial}{\partial t}u=\Delta u, \hspace{6mm}\frac{\partial}{\partial t}v=\Delta v
\end{equation}
with $|v|<u$, then we have
\begin{equation}
\n_i\n_{\bar{j}}\ln u+\frac{1}{t}g_{i\bar{j}}> \frac{\n_ih\n_{\bar{j}}h}{1-h^2},
\end{equation}
where $h=v/u$.
\end{corollary}

Taking $\epsilon=1$, we get the constrained Li-Yau-Hamilton estimate on K\"{a}hler manifolds under the K\"{a}hler-Ricci flow.
\begin{corollary}
Let $(M,g_{i\bar{j}}(t))$ be a compact solution of K\"{a}hler-Ricci flow
\begin{equation}
 \frac{\partial}{\partial t}g_{i\bar{j}}=-R_{i\bar{j}}.
 \end{equation}
 If $u$ and $v$ are the solutions of the forward conjugate heat equation
\begin{equation}
\frac{\partial}{\partial t}u=\Delta u+Ru, \hspace{6mm}\frac{\partial}{\partial t}v=\Delta v+Rv
\end{equation}
with $|v|<u$, then we have
\begin{equation}
\n_i\n_{\bar{j}}\ln u+\frac{1}{t}g_{i\bar{j}}+R_{i\bar{j}}> \frac{\n_ih\n_{\bar{j}}h}{1-h^2},
\end{equation}
where $h=v/u$.
\end{corollary}

\begin{remark}If $h=C\in(-1,1)$ is a constant, then $\nabla h=0$ and we recover the Li-Yau-Hamilton estimates in \cite{CN} and \cite{N}. Since the manifold in Theorem A is complete, we expect the result also hold on complete noncompact manifolds under suitable condition.\end{remark}

\begin{remark}
Since the constrained Li-Yau-Hamilton estimate on Riemannian manifolds helps to understand the singularities of the Ricci flow, we also expect our estimates are useful in investigating the singularities of K\"{a}hler-Ricci flow.
\end{remark}

The rest of the paper is organized as follows. We devote section 2 to the notations and basic formulas in K\"{a}hler geometry. In section 3, we prove the Theorem 1.1 by using Hamilton's maximum principle for tensors.

\section{Preliminaries}\setcounter{equation}{0}
Let $(M,g_{i\bar{j}})$ be a  K\"{a}hler manifold with K\"{a}hler metric $g_{i\bar{j}}$. The K\"{a}hler form
\begin{equation}
\omega=\frac{\sqrt{-1}}{2}g_{i\bar{j}}dz^i\wedge d\bar{z}^j
\end{equation}
is a closed real $(1,1)$-form, which is equivalent to
\begin{equation}
\frac{\partial g_{i\bar{j}}}{\partial z^k}=\frac{\partial g_{k\bar{j}}}{\partial z^i},\hspace{4mm}\frac{\partial g_{i\bar{j}}}{\partial \bar{z}^k}=\frac{\partial g_{i\bar{k}}}{\partial \bar{z}^j}.
\end{equation}

The Christoffel symbols of the metric $g_{i\bar{j}}$ is defined as
\begin{equation}
\Gamma^k_{ij}=g^{k\bar{l}}\frac{\partial g_{i\bar{l}}}{\partial z^j}, \hspace{4mm} \Gamma^{\bar{l}}_{\bar{i}\bar{j}}=g^{k\bar{l}}\frac{\partial g_{k\bar{j}}}{\partial \bar{z}^i},
\end{equation}
where $g^{i\bar{j}}=(g_{i\bar{j}})^{-1}$. It is easy to see that $\Gamma^k_{ij}$ is symmetric in $i$ and $j$ and $\Gamma^{\bar{l}}_{\bar{i}\bar{j}}$ is symmetric in $\bar{i}$ and $\bar{j}$.

The curvature tensor of the metric $g_{i\bar{j}}$ is defined by
\begin{equation}
R^j_{ik\bar{l}}=\frac{\partial\Gamma^j_{ik}}{\partial \bar{z}^l},\hspace{4mm}
R_{i\bar{j}k\bar{l}}=g_{p\bar{j}}R^p_{ik\bar{l}}.
\end{equation}
It is easy to see that $R_{i\bar{j}k\bar{l}}$ is symmetric in $i$ and $k$, in $\bar{j}$ and $\bar{l}$ and in pairs $i\bar{j}$ and $k\bar{l}$. The second Bianchi identity in K\"{a}hler case reduces to
\begin{equation}
\nabla_p R_{i\bar{j}k\bar{l}}=\nabla_k R_{i\bar{j}p\bar{l}}, \hspace{4mm}\nabla_{\bar{q}}R_{i\bar{j}k\bar{l}}=\nabla_{\bar{l}}R_{i\bar{j}p\bar{q}}.
\end{equation}

We say that $M$ has nonnegative holomorphic bisectional curvature if
\begin{equation}
R_{i\bar{j}k\bar{l}}v^iv^{\bar{j}}w^kw^{\bar{l}}\geq 0
\end{equation}
for all nonzero vectors $v$ and $w$ in the holomorphic tangent bundle $T_xM$ at $x\in M$.

The Ricci tensor of the metric $g_{i\bar{j}}$ is obtained by taking the trace of the curvature tensor:
\begin{equation}
R_{i\bar{j}}=g^{k\bar{l}}R_{{i\bar{j}k\bar{l}}}
\end{equation}
and the scalar curvature is
\begin{equation}
R=g^{i\bar{j}}R_{i\bar{j}}.
\end{equation}

Finally, we give the commutation formulas for covariant differentiations in K\"{a}hler geometry. Covariant differentiations of the same type can be commuted freely, e.g.,
\begin{equation}
\nabla_k\nabla_jv_i=\nabla_j\nabla_kv_i,\hspace{4mm}\nabla_{\bar{k}}\nabla_{\bar{j}}v_i=\nabla_{\bar{j}}\nabla_{\bar{k}}v_i.
\end{equation}
But we shall need following formulas when commuting covariant derivatives of different types:
\begin{equation}
\nabla_k\nabla_{\bar{j}}v_i=\nabla_{\bar{j}}\nabla_kv_i-R_{k{\bar{j}}i{\bar{l}}}v_l.
\end{equation}

\section{Constrained Li-Yau-Hamilton estimate}
\setcounter{equation}{0}
In this section, we will give the proof the Theorem 1.1. First we give some lemmas used in the proof. It should be pointed out that all the computations are taken in normal coordinate.
\begin{lemma} Set $L=\ln u$, then the evolution equation of $\n_{i}\n_{\bar{j}}L$ is
\begin{equation}
\begin{split}
\frac{\partial}{\partial t}\nabla_i\nabla_{\bar{j}}L=&\Delta\left(\nabla_i\nabla_{\bar{j}}L\right)+R_{l\bar{j}i\bar{k}}\nabla_k\nabla_{\bar{l}}L
+R_{i\bar{j}k\bar{l}}{\nabla_{l}L\nabla_{\bar{k}}L}\\&+\nabla_{k}L\nabla_{\bar{k}}\nabla_i\nabla_{\bar{j}}L+\nabla_{\bar{k}}L\nabla_k\nabla_i\nabla_{\bar{j}}L
\\&+\nabla_{\bar{j}}\nabla_{k}L\nabla_i\nabla_{\bar{k}}L+\nabla_i\nabla_{k}L\nabla_{\bar{j}}\nabla_{\bar{k}}L
\\&-\frac{1}{2}\left(R_{l\bar{j}}\nabla_i\nabla_{\bar{l}}L+R_{i\bar{l}}\nabla_l\nabla_{\bar{j}}L\right)\\&+\epsilon\left(\Delta R_{i{\bar{j}}}+R_{i{\bar{j}}k{\bar{l}}}{ R_{{\bar{k}}l}}-R_{i{\bar{p}}}R_{p{\bar{j}}}\right).
\end{split}
\end{equation}
\end{lemma}
Proof: By
\begin{equation}
\frac{\partial}{\partial t}L=\frac{\partial}{\partial t}\ln u=\frac{u_t}{u}=\frac{1}{u}\left(\Delta u+\epsilon Ru\right)=\frac{\Delta u}{u}+\epsilon R,
\end{equation}
and
\begin{equation}
\Delta L=\n_i\n_{\bar{i}}L=\n_i\frac{\n_{\bar{i}}u}{u}=\frac{\Delta u}{u}-\frac{|\nabla u|^2}{u^2},
\end{equation}
we get the evolution equation of $L$
\begin{equation}
\frac{\partial}{\partial t}L=\Delta L+\frac{|\nabla u|^2}{u^2}+\epsilon R=\Delta L+|\n L|^2+\epsilon R.
\end{equation}
We claim that if a function $A$ satisfies the equation
\begin{equation}
\frac{\partial}{\partial t}A=\Delta {A}+B,
\end{equation}
then $\nabla_i\nabla_{\bar{j}}A$ satisfies
\begin{equation}
\begin{split}
\frac{\partial}{\partial t}\left(\nabla_i\nabla_{\bar{j}}A\right)=&\Delta\left(\nabla_i\nabla_{\bar{j}}A\right)+R_{i\bar{j}l\bar{k}}\nabla_k\nabla_{\bar{l}}A\\&-\frac{1}{2}\left(R_{l\bar{j}}\nabla_i\nabla_{\bar{l}}A+R_{i\bar{l}}\nabla_l\nabla_{\bar{j}}A\right)
+\nabla_i\nabla_{\bar{j}}B.
\end{split}
\end{equation}
In fact a straightforward computation leads to
\begin{equation}
\begin{split}
\frac{\partial}{\partial t}\left(\nabla_i\nabla_{\bar{j}}A\right)&=\nabla_i\nabla_{\bar{j}}\left(\frac{\partial}{\partial t}A\right)\\
&=\nabla_i\nabla_{\bar{j}}\left(\Delta {A}+B\right)\\
&=\nabla_i\nabla_{\bar{j}}\left(\Delta {A}\right)+\nabla_i\nabla_{\bar{j}}B. \label{evol1}
\end{split}
\end{equation}
It is sufficient to compute
 $$\nabla_i\nabla_{\bar{j}}\left(\Delta {A}\right)=\frac{1}{2}\nabla_i\nabla_{\bar{j}}\left(\nabla_k\nabla_{\bar{k}}A\right)+\nabla_i\nabla_{\bar{j}}\left(\nabla_{\bar{k}}\nabla_{k}A\right).$$
 By using communication formula, we get
\begin{equation}
\begin{split}
\nabla_i\nabla_{\bar{j}}\left(\nabla_k\nabla_{\bar{k}}A\right)=&\nabla_i\left(\nabla_k\nabla_{\bar{j}}\nabla_{\bar{k}}A-R_{l\bar{j}}\nabla_{\bar{l}}A\right)\\
=&\nabla_k\nabla_i\nabla_{\bar{k}}\nabla_{\bar{j}}A-\nabla_i{R_{l\bar{j}}}\nabla_{\bar{l}}A-R_{l\bar{j}}\nabla_i\nabla_{\bar{l}}A\\
=&\nabla_k\left(\nabla_{\bar{k}}\nabla_i\nabla_{\bar{j}}A+R_{i\bar{k}l\bar{j}}\nabla_{\bar{l}}A\right)-\nabla_i{R_{l\bar{j}}}\nabla_{\bar{l}}A-R_{l\bar{j}}\nabla_i\nabla_{\bar{l}}A\\
=&\nabla_k\nabla_{\bar{k}}\nabla_i\nabla_{\bar{j}}A+R_{i\bar{k}l\bar{j}}\nabla_k\nabla_{\bar{l}}A+\nabla_k{R_{i\bar{k}l\bar{j}}}\nabla_{\bar{l}}A
-\nabla_i{R_{l\bar{j}}}\nabla_{\bar{l}}A-R_{l\bar{j}}\nabla_i\nabla_{\bar{l}}A\\
=&\nabla_k\nabla_{\bar{k}}\nabla_i\nabla_{\bar{j}}A+R_{i\bar{j}l\bar{k}}\nabla_k\nabla_{\bar{l}}A-R_{l\bar{j}}\nabla_i\nabla_{\bar{l}}A  \label{evol2}
\end{split}
\end{equation}
and
\begin{equation}
\begin{split}
\nabla_i\nabla_{\bar{j}}\left(\nabla_{\bar{k}}\nabla_{k}A\right)=&\nabla_i\left(\nabla_{\bar{k}}\nabla_{\bar{j}}\nabla_{k}A\right)\\=&\nabla_i\left(\nabla_{\bar{k}}\nabla_{k}\nabla_{\bar{j}}A\right)\\
=&\nabla_{\bar{k}}\nabla_i\nabla_k\nabla_{\bar{j}}A+R_{i\bar{k}l\bar{j}}\nabla_k\nabla_{\bar{l}}A-R_{i\bar{k}k\bar{l}}\nabla_l\nabla_{\bar{j}}A\\
=&\nabla_{\bar{k}}\nabla_k\left(\nabla_i\nabla_{\bar{j}}A\right)+R_{i\bar{j}l\bar{k}}\nabla_k\nabla_{\bar{l}}A-R_{i\bar{l}}\nabla_l\nabla_{\bar{j}}A. \label{evol3}
\end{split}
\end{equation}
Inserting (\ref{evol2}) and (\ref{evol3}) into (\ref{evol1}), we get the evolution equation of $\nabla_i\nabla_{\bar{j}}A$.

Applying this formula to $A=L=\ln u$, we obtain
\begin{equation}
\begin{split}
\frac{\partial}{\partial t}\left(\nabla_i\nabla_{\bar{j}}L\right)=&\Delta\left(\nabla_i\nabla_{\bar{j}}L\right)+R_{i\bar{j}l\bar{k}}\nabla_k\nabla_{\bar{l}}L\\
&-\frac{1}{2}\left(R_{l\bar{j}}\nabla_i\nabla_{\bar{l}}L+R_{i\bar{l}}\nabla_l\nabla_{\bar{j}}L\right)+\nabla_i\nabla_{\bar{j}}|\nabla L|^2
+\epsilon \nabla_i\nabla_{\bar{j}}R.
\end{split}
\end{equation}
By Lemma 1.1 in \cite{N}, we get
\begin{equation}
 \nabla_i\nabla_{\bar{j}} R=\Delta R_{i{\bar{j}}}+R_{i{\bar{j}}k{\bar{l}}}R_{{\bar{k}}l}-R_{i{\bar{p}}}R_{p{\bar{j}}}.
\end{equation}
Now the lemma follows from the following computation
\begin{equation}
\begin{split}
\nabla_i\nabla_{\bar{j}}|\nabla L|^2=&\nabla_i\nabla_{\bar{j}}\left(\nabla_{k}L\nabla_{\bar{k}}L\right)\\=&\nabla_i\left(\nabla_{\bar{j}}\nabla_{k}L\nabla_{\bar{k}}L+\nabla_{k}L\nabla_{\bar{j}}\nabla_{\bar{k}}L\right)
\\=&\nabla_i\nabla_{\bar{j}}\nabla_{k}L\nabla_{\bar{k}}L+\nabla_{\bar{j}}\nabla_{k}L\nabla_i\nabla_{\bar{k}}L+\nabla_i\nabla_{k}L\nabla_{\bar{j}}\nabla_{\bar{k}}L
+\nabla_{k}L\nabla_i\nabla_{\bar{j}}\nabla_{\bar{k}}L
\\=&\nabla_{\bar{k}}L\nabla_k\left(\nabla_i\nabla_{\bar{j}}L\right)+\nabla_{\bar{j}}\nabla_{k}L\nabla_i\nabla_{\bar{k}}L\\&+\nabla_i\nabla_{k}L\nabla_{\bar{j}}\nabla_{\bar{k}}L
+\nabla_{k}L\nabla_{\bar{k}}\left(\nabla_i\nabla_{\bar{j}}L\right)+R_{i\bar{j}k\bar{l}}\nabla_{l}L\nabla_{\bar{k}}L.
\end{split}
\end{equation}

\begin{lemma}(see \cite{N}) If K\"{a}hler metrics $g_{i\bar{j}}(t)$ satisfy the $\epsilon$-K\"{a}hler-Ricci flow, then the Ricci tensor satisfies
\begin{equation}
\frac{\partial}{\partial t}\left(\epsilon R_{i\bar{j}}\right)=\epsilon^{2}\left(\Delta R_{i{\bar{j}}}+R_{i{\bar{j}}k{\bar{l}}}R_{{\bar{k}}l}-R_{i{\bar{p}}}R_{p{\bar{j}}}\right).
\end{equation}
\end{lemma}

\begin{lemma}Let $h=v/u$, then
\begin{equation}
\begin{split}
\frac{\partial}{\partial t}\left(\frac{\n_{i}h\n_{\bar{j}}h}{1-h^2}\right)
=&\Delta\left(\frac{\n_{i}h\n_{\bar{j}}h}{1-h^2}\right)+\n_{k}L\n_{\bar{k}}\left(\frac{\n_{i}h\n_{\bar{j}}h}{1-h^2}\right)+\n_{\bar{k}}L\n_{k}\left(\frac{\n_{i}h\n_{\bar{j}}h}{1-h^2}\right)\\
&-\frac{1}{1-h^2}\left(\n_{i}\n_{k}h+\frac{2h\n_{i}h\n_{k}h}{1-h^2}\right)\left(\n_{\bar{j}}\n_{\bar{k}}h+\frac{2h\n_{\bar{j}}h\n_{\bar{k}}h}{1-h^2}\right)\\
&-\frac{1}{1-h^2}\left(\n_{i}\n_{\bar{k}}h+\frac{2h\n_{i}h\n_{\bar{k}}h}{1-h^2}\right)\left(\n_{\bar{j}}\n_{k}h+\frac{2h\n_{\bar{j}}h\n_{k}h}{1-h^2}\right)\\
&+\n_{i}\n_{k}L\left(\frac{\n_{\bar{j}}h\n_{\bar{k}}h}{1-h^2}\right)+\n_{i}\n_{\bar{k}}L\left(\frac{\n_{\bar{j}}h\n_{k}h}{1-h^2}\right)\\&+\n_{\bar{j}}\n_{k}L\left(\frac{\n_{i}h\n_{\bar{k}}h}{1-h^2}\right)
+\n_{\bar{j}}\n_{\bar{k}}L\left(\frac{\n_{i}h\n_{k}h}{1-h^2}\right)
\\&-\frac{1}{2}R_{i\bar{l}}\left(\frac{\n_{\bar{j}}h\n_{l}h}{1-h^2}\right)-\frac{1}{2}R_{k\bar{j}}\left(\frac{\n_{\bar{k}}h\n_{i}h}{1-h^2}\right)\\&-\frac{2\n_{i}h\n_{\bar{j}}h}{\left(1-h^2\right)^2}|\n h|^2.
\end{split}
\end{equation}
\end{lemma}
Proof: We divide our proof in fours steps. First, we need to give the evolution equation of $h$. By the definition of $h$, we have
\begin{equation}
\frac{\partial}{\partial t}h=\frac{\partial}{\partial t}\left(\frac{v}{u}\right)=\frac{v_t}{u}-\frac{vu_t}{u^2}=\frac{1}{u^2}\left(u\Delta v-v\Delta u\right)\label{i1}
\end{equation}
and
\begin{equation}
\begin{split}
\Delta h=\nabla_k\nabla_{\bar{k}}h=&\nabla_k\frac{u\nabla_{\bar{k}}v-v\nabla_{\bar{k}}u}{u^2}
\\=&\frac{\nabla_ku\nabla_{\bar{k}}v+u\nabla_k\nabla_{\bar{k}}v-\nabla_kv\nabla_{\bar{k}}u-v\nabla_k\nabla_{\bar{k}}u}{u^2}
\\&-\frac{2\nabla_ku\left(u\nabla_{\bar{k}}v-v\nabla_{\bar{k}}u\right)}{u^3}
\\=&\frac{u\Delta v-v\Delta u}{u^2}-\frac{\nabla_{k}u\nabla_{\bar{k}}v+\nabla_kv\nabla_{\bar{k}}u}{u^2}+\frac{2v\nabla_{\bar{k}}u\nabla_{k}u}{u^3}
\\=&\frac{u\Delta v-v\Delta u}{u^2}-\nabla_{k}L\nabla_{\bar{k}}h-\nabla_{\bar{k}}L\nabla_{k}h.\label{i2}
\end{split}
\end{equation}
Combining (\ref{i1}) and (\ref{i2}), we obtain the evolution of $h$
\begin{equation}
\begin{split}
\frac{\partial}{\partial t}h=&\Delta h+\nabla_{k}L\nabla_{\bar{k}}h+\nabla_{\bar{k}}L\nabla_{k}h.
\end{split}
\end{equation}

The next thing to do in the proof is computing the evolution equation of $1-h^2$. By
\begin{equation}
\begin{split}
\frac{\partial}{\partial t}\left(1-h^2\right)=&-2h\frac{\partial}{\partial t}h\\=&-2h\left(\Delta h+\nabla_{k}L\nabla_{\bar{k}}h+\nabla_{\bar{k}}L\nabla_{k}h\right)
\end{split}
\end{equation}
and
\begin{equation}
\Delta\left(1-h^2\right)=-2|\nabla h|^2-2h\Delta h,
\end{equation}
we have that
\begin{equation}
\begin{split}
\frac{\partial}{\partial t}\left(1-h^2\right)=&\Delta\left(1-h^2\right)+2|\nabla h|^2{}-2h\left(\nabla_{k}L\nabla_{\bar{k}}h+\nabla_{\bar{k}}L\nabla_{k}h\right)\\=&\Delta\left(1-h^2\right)+2|\nabla h|^2+\nabla_{k}L\nabla_{\bar{k}}\left(1-h^2\right)+\nabla_{\bar{k}}L\nabla_{k}\left(1-h^2\right). \label{gradient evol1}
\end{split}
\end{equation}

Another step in the proof is computing the evolution equation of $\nabla_{i}h\nabla_{\bar{j}}h$. Take partial derivative with respective to $t$
\begin{equation}
\frac{\partial}{\partial t}\left(\nabla_{i}h\nabla_{\bar{j}}h\right)=\nabla_{i}\left(\frac{\partial}{\partial t}h\right)\nabla_{\bar{j}}h+\nabla_{i}h\nabla_{\bar{j}}\left(\frac{\partial}{\partial t}h\right).
\end{equation}
Direct computation yields
\begin{equation}
\begin{split}
\nabla_{i}\left(\frac{\partial}{\partial t}h\right)\nabla_{\bar{j}}h=&\nabla_{i}\left(\Delta h+\nabla_{k}L\nabla_{\bar{k}}h+\nabla_{\bar{k}}L\nabla_{k}h\right)\nabla_{\bar{j}}h
\\=&\nabla_{\bar{j}}h\Delta\left(\nabla_{i}h\right)-\frac{1}{2}R_{i\bar{l}}\nabla_{l}h\nabla_{\bar{j}}h+\nabla_i\nabla_{k}L\nabla_{\bar{k}}h\nabla_{\bar{j}}h
\\&+\nabla_{k}L\nabla_i\nabla_{\bar{k}}h\nabla_{\bar{j}}h+\nabla_i\nabla_{\bar{k}}L\nabla_{k}h\nabla_{\bar{j}}h+\nabla_{\bar{k}}L\nabla_i\nabla_{k}h\nabla_{\bar{j}}h
\end{split}\label{i3}
\end{equation}
and
\begin{equation}
\begin{split}
\nabla_{\bar{j}}\left(\frac{\partial}{\partial t}h\right)\nabla_{i}h=&\nabla_{\bar{j}}\left(\Delta h+\nabla_{k}L\nabla_{\bar{k}}h+\nabla_{\bar{k}}L\nabla_{k}h\right)\nabla_{i}h
\\=&\n_{\bar{j}}\Delta h\n_{i}h+\n_{\bar{j}}\n_{k}L\n_{\bar{k}}h\n_{i}h+\n_{k}L\n_{\bar{j}}\n_{\bar{k}}h\n_{i}h\\
&+\n_{\bar{j}}\n_{\bar{k}}L\n_{k}h\n_{i}h+\n_{\bar{k}}L\n_{\bar{j}}\n_{k}h\n_{i}h\\
=&\Delta\left(\n_{\bar{j}}h\right)\n_{i}h-\frac{1}{2}R_{k\bar{j}}\n_{\bar{k}}h\n_{i}h\\
&+\n_{\bar{j}}\n_{k}L\n_{\bar{k}}h\n_{i}h+\n_{k}L\n_{\bar{j}}\n_{\bar{k}}h\n_{i}h\\
&+\n_{\bar{j}}\n_{\bar{k}}L\n_{k}h\n_{i}h+\n_{\bar{k}}L\n_{\bar{j}}\n_{k}h\n_{i}h.
\end{split}\label{i4}
\end{equation}
Combining (\ref{i3}) and (\ref{i4}), we have that
\begin{equation}
\begin{split}
\frac{\partial}{\partial t}\left(\n_{i}h\n_{\bar{j}}h\right)=&\Delta\left(\n_{i}h\right)\n_{\bar{j}}h+\Delta\left(\n_{\bar{j}}h\right)\n_{i}h-\frac{1}{2}R_{i\bar{l}}\n_{l}h\n_{\bar{j}}h\\
&-\frac{1}{2}R_{k\bar{j}}\n_{\bar{k}}h\n_{i}h+\n_{i}\n_{k}L\n_{\bar{j}}h\n_{\bar{k}}h+\n_{i}\n_{\bar{k}}L\n_{k}h\n_{\bar{j}}h\\
&+\n_{k}L\n_{\bar{k}}\left(\n_{i}h\n_{\bar{j}}h\right)+\n_{\bar{k}}L\n_{k}\left(\n_{i}h\n_{\bar{j}}h\right)
\\&+\n_{\bar{j}}\n_{k}L\n_{\bar{k}}h\n_{i}h+\n_{\bar{j}}\n_{\bar{k}}L\n_{k}h\n_{i}h.
\end{split}
\end{equation}
The Laplacian of $\nabla_{i}h\nabla_{\bar{j}}h$ is
\begin{equation}
\begin{split}
\Delta\left(\n_{i}h\n_{\bar{j}}h\right)
=&\frac{1}{2}\left(\n_{k}\n_{\bar{k}}+\n_{\bar{k}}\n_{k}\right)\left(\n_{i}h\n_{\bar{j}}h\right)\\
=&\frac{1}{2}\n_{k}\left(\n_{\bar{k}}\n_{i}h\n_{\bar{j}}h+\n_{i}h\n_{\bar{k}}\n_{\bar{j}}h\right)
\\&+\frac{1}{2}\n_{\bar{k}}\left(\n_{k}\n_{i}h\n_{\bar{j}}h+\n_{i}h\n_{k}\n_{\bar{j}}h\right)\\
=&\Delta\left(\n_{i}h\right)\n_{\bar{j}}h+\n_{i}h\Delta\left(\n_{\bar{j}}h\right)\\&+\n_{k}\n_{i}h\n_{\bar{k}}\n_{\bar{j}}h+\n_{\bar{k}}\n_{i}h\n_{k}\n_{\bar{j}}h.
\end{split}
\end{equation}
Collecting the two equalities above gives that
\begin{equation}
\begin{split}
\frac{\partial}{\partial t}\left(\n_{i}h\n_{\bar{j}}h\right)=&\Delta\left(\n_{i}h\n_{\bar{j}}h\right)-\n_{k}\n_{i}h\n_{\bar{k}}\n_{\bar{j}}h-\n_{\bar{k}}\n_{i}h\n_{k}\n_{\bar{j}}h
\\&+\n_{i}\n_{k}L\n_{\bar{j}}h\n_{\bar{k}}h
+\n_{i}\n_{\bar{k}}L\n_{k}h\n_{\bar{j}}h\\&+\n_{\bar{j}}\n_{k}L\n_{\bar{k}}h\n_{i}h+\n_{\bar{j}}\n_{\bar{k}}L\n_{k}h\n_{i}h\\
&+\n_{k}L\n_{\bar{k}}\left(\n_{i}h\n_{\bar{j}}h\right)+\n_{\bar{k}}L\n_{k}\left(\n_{i}h\n_{\bar{j}}h\right)\\&
-\frac{1}{2}R_{i\bar{l}}\n_{\bar{j}}h\n_{l}h-\frac{1}{2}R_{k\bar{j}}\n_{\bar{k}}h\n_{i}h. \label{gradient evol2}
\end{split}
\end{equation}

Finally, we give the evolution of $\frac{\n_{i}h\n_{\bar{j}}h}{1-h^2}$. With the help of (\ref{gradient evol1}) and (\ref{gradient evol2}), we have that
\begin{equation}
\begin{split}
\frac{\partial}{\partial t}\left(\frac{\n_{i}h\n_{\bar{j}}h}{1-h^2}\right)=&\frac{1}{1-h^2}\frac{\partial}{\partial t}{\left(\n_{i}h\n_{\bar{j}}h\right)}-\frac{\n_{i}h\n_{\bar{j}}h}{\left(1-h^2\right)^2}\frac{\partial}{\partial t}\left(1-h^2\right){}\\=&\frac{1}{1-h^2}[\Delta\left(\n_{i}h\n_{\bar{j}}h\right)-\n_{k}\n_{i}h\n_{\bar{k}}\n_{\bar{j}}h-\n_{\bar{k}}\n_{i}h\n_{k}\n_{\bar{j}}h{}\\&+\n_{i}\n_{k}L\n_{\bar{j}}h\n_{\bar{k}}h
+\n_{i}\n_{\bar{k}}L\n_{k}h\n_{\bar{j}}h+\n_{\bar{j}}\n_{k}L\n_{\bar{k}}h\n_{i}h\\&+\n_{\bar{j}}\n_{\bar{k}}L\n_{i}h\n_{k}h+\n_{k}L\n_{\bar{k}}(\n_{i}h\n_{\bar{j}}h){}\\
&+\n_{\bar{k}}L\n_{k}\left(\n_{i}h\n_{\bar{j}}h\right)-\frac{1}{2}R_{i\bar{l}}\n_{\bar{j}}h\n_{l}h-\frac{1}{2}R_{k\bar{j}}\n_{\bar{k}}h\n_{i}h]{}\\
&-\frac{\n_{i}h\n_{\bar{j}}h}{\left(1-h^2\right)^2}[\Delta\left(1-h^2\right)+2|\n h|^2+\n_{k}L\n_{\bar{k}}\left(1-h^2\right) \\&+\n_{\bar{k}}L\n_{k}\left(1-h^2\right)].
\end{split}
\end{equation}
Straightforward computation yields
\begin{equation}
\begin{split}
\Delta\left(\frac{\n_{i}h\n_{\bar{j}}h}{1-h^2}\right)=&\frac{\Delta\left(\n_{i}h\n_{\bar{j}}h\right)}{1-h^2}-\frac{\n_{i}h\n_{\bar{j}}h}{\left(1-h^2\right)^2}\Delta\left(1-h^2\right)-\frac{\n_{k}\left(\n_{i}h\n_{\bar{j}}h\right)\n_{\bar{k}}\left(1-h^2\right)}{\left(1-h^2\right)^2}\\&
-\frac{\n_{\bar{k}}\left(\n_{i}h\n_{\bar{j}}h\right)\n_{k}\left(1-h^2\right)}{\left(1-h^2\right)^2}+\frac{2\n_{i}h\n_{\bar{j}}h[\n_{k}\left(1-h^2\right)\n_{\bar{k}}\left(1-h^2\right)]}{\left(1-h^2\right)^3}\\
=&\frac{\Delta\left(\n_{i}h\n_{\bar{j}}h\right)}{1-h^2}-\frac{\n_{i}h\n_{\bar{j}}h}{\left(1-h^2\right)^2}\Delta\left(1-h^2\right)\\
&+\frac{2h\n_{k}\left(\n_{i}h\n_{\bar{j}}h\right)\n_{\bar{k}}h+2h\n_{\bar{k}}\left(\n_{i}h\n_{\bar{j}}h\right)\n_{k}h}{\left(1-h^2\right)^2}\\&+\frac{8h^2\n_{i}h\n_{\bar{j}}h|\n h|^2}{\left(1-h^2\right)^3},
\end{split}
\end{equation}
so we have that
\begin{equation}
\begin{split}
\frac{\partial}{\partial t}\left(\frac{\n_{i}h\n_{\bar{j}}h}{1-h^2}\right)
=&\Delta\left(\frac{\n_{i}h\n_{\bar{j}}h}{1-h^2}\right)+\frac{1}{1-h^2}(-\n_{k}\n_{i}h\n_{\bar{k}}\n_{\bar{j}}h\\&-\n_{\bar{k}}\n_{i}h\n_{k}\n_{\bar{j}}h+\n_{i}\n_{k}L\n_{\bar{j}}h\n_{\bar{k}}h\\
&+\n_{i}\n_{\bar{k}}L\n_{\bar{j}}h\n_{k}h+\n_{\bar{j}}\n_{k}L\n_{\bar{k}}h\n_{i}h\\
&+\n_{\bar{j}}\n_{\bar{k}}L\n_{i}h\n_{k}h-\frac{1}{2}R_{i\bar{l}}\n_{\bar{j}}h\n_{l}h-\frac{1}{2}R_{k\bar{j}}\n_{\bar{k}}h\n_{i}h)
\\&-\frac{2\n_{i}h\n_{\bar{j}}h}{\left(1-h^2\right)^2}|\n h|^2+\n_{k}L\n_{\bar{k}}\left(\frac{\n_{i}h\n_{\bar{j}}h}{1-h^2}\right)\\
&+\n_{\bar{k}}L\n_{k}\left(\frac{\n_{i}h\n_{\bar{j}}h}{1-h^2}\right)-\frac{8h^2\n_{i}h\n_{\bar{j}}h|\n h|^2}{\left(1-h^2\right)^3}\\
&-\frac{2h\n_{k}\left(\n_{i}h\n_{\bar{j}}h\right)\n_{\bar{k}}h+2h\n_{\bar{k}}\left(\n_{i}h\n_{\bar{j}}h\right)\n_{k}h}{\left(1-h^2\right)^2}.
\end{split}
\end{equation}
Rearranging terms leads to
\begin{equation}
\begin{split}
\frac{\partial}{\partial t}\left(\frac{\n_{i}h\n_{\bar{j}}h}{1-h^2}\right)
=&\Delta\left(\frac{\n_{i}h\n_{\bar{j}}h}{1-h^2}\right)+\n_{k}L\n_{\bar{k}}\left(\frac{\n_{i}h\n_{\bar{j}}h}{1-h^2}\right)+\n_{\bar{k}}L\n_{k}\left(\frac{\n_{i}h\n_{\bar{j}}h}{1-h^2}\right)\\
&-\frac{1}{1-h^2}\left(\n_{i}\n_{k}h+\frac{2h\n_{i}h\n_{k}h}{1-h^2}\right)\left(\n_{\bar{j}}\n_{\bar{k}}h+\frac{2h\n_{\bar{j}}h\n_{\bar{k}}h}{1-h^2}\right)\\
&-\frac{1}{1-h^2}\left(\n_{i}\n_{\bar{k}}h+\frac{2h\n_{i}h\n_{\bar{k}}h}{1-h^2}\right)\left(\n_{\bar{j}}\n_{k}h+\frac{2h\n_{\bar{j}}h\n_{k}h}{1-h^2}\right)\\
&+\frac{1}{1-h^2}(\n_{i}\n_{k}L\n_{\bar{j}}h\n_{\bar{k}}h+\n_{i}\n_{\bar{k}}L\n_{\bar{j}}h\n_{k}h+\n_{\bar{j}}\n_{k}L\n_{i}h\n_{\bar{k}}h\\&+\n_{\bar{j}}\n_{\bar{k}}L\n_{i}h\n_{k}h
-\frac{1}{2}R_{i\bar{l}}\n_{\bar{j}}h\n_{l}h-\frac{1}{2}R_{k\bar{j}}\n_{\bar{k}}h\n_{i}h)\\
&-\frac{2\n_{i}h\n_{\bar{j}}h}{\left(1-h^2\right)^2}|\n h|^2.
\end{split}
\end{equation}
This complete the proof of lemma 3.

We also need following result, which is an easy consequence of Cao's differential Harnack inequality for the K\"{a}hler-Ricci flow \cite{C2}.
\begin{lemma} Let $(M,g_{i\bar{j}}(t))$ be a solution of $\epsilon$-K\"{a}hler-Ricci flow with nonnegative bisectional curvature and $u$ be a solution of (9), then there holds
\begin{eqnarray*}
Y_{i{\bar{j}}}&=&\Delta R_{i{\bar{j}}}+R_{i{\bar{j}}{k}{\bar{l}}}R_{{\bar{k}}l}-\left(\frac{\n_{\bar{k}}L}{\epsilon}\n_{k}R_{i{\bar{j}}}+\frac{\n_{k}L}{\epsilon}\n_{\bar{k}}R_{i{\bar{j}}}\right)
+R_{i{\bar{j}}l{\bar{k}}}\frac{\n_{k}L}{\epsilon}\frac{\n_{\bar{l}}L}{\epsilon}+\frac{R_{i{\bar{j}}}}{\epsilon t}\geq 0.
\end{eqnarray*}
\end{lemma}

With the help of the preceding four lemmas we are now in the position to give the proof of the theorem.\newline

{\bf Proof of the Theorem 1.1}
Let \begin{equation}
P_{i{\bar{j}}}=\n_{i}\n_{\bar{j}}L+\epsilon R_{i\bar{j}}-\frac{\n_{i}h\n_{\bar{j}}h}{1-h^2}.
\end{equation}
Combine the lemmas together to get
\begin{equation}
\begin{split}
\frac{\partial}{\partial t}P_{i{\bar{j}}}=&\Delta P_{i{\bar{j}}} +\nabla_{k}L\nabla_{\bar{k}}P_{i{\bar{j}}}+\nabla_{\bar{k}}L\nabla_{k}P_{i{\bar{j}}}+R_{i{\bar{j}}k{\bar{l}}}P_{l{\bar{k}}}+\epsilon^2Y_{i{\bar{j}}}-\epsilon^2R_{i\bar{p}}R_{p\bar{j}}
\\&-\frac{1}{2}\left(R_{l{\bar{j}}}P_{i{\bar{l}}}+R_{i{\bar{l}}}P_{l{\bar{j}}}\right)+\frac{1}{2}R_{l{\bar{j}}}\left(\epsilon R_{i{\bar{l}}}\right)+\frac{1}{2}R_{i{\bar{l}}}\left(\epsilon R_{l{\bar{j}}}\right)-\epsilon R_{i\bar{p}}R_{p\bar{j}}\\&+\left(\nabla_{\bar{j}}\nabla_{k}L-\frac{\n_{\bar{j}}h\n_{k}h}{1-h^2}\right)\left(\nabla_{i}\nabla_{\bar{k}}L-\frac{\n_{i}h\n_{\bar{k}}h}{1-h^2}\right)
\\&+\left(\n_{i}\n_{k}L-\frac{\n_{i}h\n_{k}h}{1-h^2}\right)\left(\n_{\bar{j}}\n_{\bar{k}}L-\frac{\n_{\bar{j}}h\n_{\bar{k}}h}{1-h^2}\right)\\&+\frac{1}{1-h^2}\left(\n_{i}\n_{k}h+\frac{2h\n_{i}h\n_{k}h}{1-h^2}\right)\left(\n_{\bar{j}}\n_{\bar{k}}h+\frac{2h\n_{\bar{j}}h\n_{\bar{k}}h}{1-h^2}\right)\\
&+\frac{1}{1-h^2}\left(\n_{i}\n_{\bar{k}}h+\frac{2h\n_{i}h\n_{\bar{k}}h}{1-h^2}\right)\left(\n_{\bar{j}}\n_{k}h+\frac{2h\n_{\bar{j}}h\n_{k}h}{1-h^2}\right).
\end{split}
\end{equation}
Since the three matrices in the last three terms are nonnegative definite, we have that
\begin{equation}
\begin{split}
\frac{\partial}{\partial t}P_{i{\bar{j}}}\geq&\Delta P_{i{\bar{j}}}+\nabla_{\bar{k}}L\nabla_{k}P_{i{\bar{j}}}+\nabla_{k}L\nabla_{\bar{k}}P_{i{\bar{j}}}+R_{i{\bar{j}}k{\bar{l}}}P_{l{\bar{k}}}
\\&+\left(\nabla_{\bar{j}}\nabla_{k}L-\frac{\n_{\bar{j}}h\n_{k}h}{1-h^2}\right)\left(\nabla_{i}\nabla_{\bar{k}}L-\frac{\n_{i}h\n_{\bar{k}}h}{1-h^2}\right)
\\&+\epsilon^2Y_{i{\bar{j}}}-\epsilon^2R_{i\bar{p}}R_{p\bar{j}}-\frac{\epsilon R_{i{\bar{j}}}}{ t}-\frac{1}{2}\left(R_{l{\bar{j}}}P_{i{\bar{l}}}+R_{i{\bar{l}}}P_{l{\bar{j}}}\right).
\end{split}
\end{equation}
Therefore, by Lemma 4, we get that
\begin{equation}
\begin{split}
\frac{\partial}{\partial t}\left(P_{i{\bar{j}}}+\frac{g_{i{\bar{j}}}}{t}\right)
=&\frac{\partial}{\partial t}P_{i{\bar{j}}}+\frac{1}{t}\frac{\partial}{\partial t}g_{i{\bar{j}}}-\frac{1}{t^2}g_{i{\bar{j}}}\\\geq&\Delta\left(P_{i{\bar{j}}}+\frac{g_{i{\bar{j}}}}{t}\right)+\n_{\bar{k}}L\n_{k}\left(P_{i{\bar{j}}}+\frac{g_{i{\bar{j}}}}{t}\right)
\\
&+\n_{k}L\n_{\bar{k}}\left(P_{i{\bar{j}}}+\frac{g_{i{\bar{j}}}}{t}\right)+R_{i{\bar{j}}k{\bar{l}}}\left(P_{l{\bar{k}}}+\frac{g_{l{\bar{k}}}}{t}\right)\\&-R_{i{\bar{j}}k{\bar{l}}}\frac{g_{l{\bar{k}}}}{t}
-\frac{1}{2}\left[R_{l{\bar{j}}}\left(P_{i{\bar{l}}}+\frac{g_{i{\bar{l}}}}{t}\right)+R_{i{\bar{l}}}\left(P_{l{\bar{j}}}+\frac{g_{l{\bar{j}}}}{t}\right)\right]
\\&+\left(\nabla_{\bar{j}}\nabla_{k}L-\frac{\n_{\bar{j}}h\n_{k}h}{1-h^2}\right)\left(\nabla_{i}\nabla_{\bar{k}}L-\frac{\n_{i}h\n_{\bar{k}}h}{1-h^2}\right)\\
&-\epsilon^2R_{i\bar{p}}R_{p\bar{j}}-\frac{2\epsilon R_{i{\bar{j}}}}{ t}-\frac{1}{t^2}g_{i{\bar{j}}}.
\end{split}
\end{equation}
Let
\begin{equation}
\begin{split}
Q_{i{\bar{j}}}=P_{i{\bar{j}}}+\frac{g_{i{\bar{j}}}}{t},
\end{split}
\end{equation}
then $Q_{i{\bar{j}}}$ satisfies that
\begin{equation}
\begin{split}
\frac{\partial}{\partial t}Q_{i{\bar{j}}}\geq&\Delta Q_{i{\bar{j}}}+\left<\n L,\n Q_{i{\bar{j}}}\right>+\left<\n Q_{i{\bar{j}}},\n L\right>+R_{i{\bar{j}}k{\bar{l}}}Q_{l{\bar{k}}}\\&-\frac{1}{2}\left(R_{l{\bar{j}}}Q_{i{\bar{l}}}+R_{i{\bar{l}}}Q_{l{\bar{j}}}\right)
+P_{i{\bar{k}}}P_{k{\bar{j}}}-\epsilon R_{k{\bar{j}}}P_{i{\bar{k}}}\\&-\epsilon R_{i{\bar{k}}}\left(P_{k{\bar{j}}}-\epsilon R_{k{\bar{j}}}\right)-\epsilon^2R_{i\bar{p}}R_{p\bar{j}}-\frac{2\epsilon R_{i{\bar{j}}}}{ t}-\frac{1}{t^2}g_{i{\bar{j}}}\\
=&\Delta Q_{i{\bar{j}}}+\left<\n L,\n Q_{i{\bar{j}}}\right>+\left<\n Q_{i{\bar{j}}},\n L\right>+R_{i{\bar{j}}k{\bar{l}}}Q_{l{\bar{k}}}\\&-\left(\frac{1}{2}+\epsilon\right)\left(R_{l{\bar{j}}}Q_{i{\bar{l}}}+R_{i{\bar{l}}}Q_{l{\bar{j}}}\right)
+\left(P_{k{\bar{j}}}+\frac{g_{k{\bar{j}}}}{t}\right)\left(P_{i{\bar{k}}}-\frac{g_{i{\bar{k}}}}{t}\right)
\\=&\Delta Q_{i{\bar{j}}}+\left<\n L,\n Q_{i{\bar{j}}}\right>+\left<\n Q_{i{\bar{j}}},\n L\right>+R_{i{\bar{j}}k{\bar{l}}}Q_{l{\bar{k}}}\\&-\left(\frac{1}{2}+\epsilon\right)\left(R_{p{\bar{j}}}Q_{i{\bar{p}}}+R_{i{\bar{p}}}Q_{p{\bar{j}}}\right)
+\left(Q_{i{\bar{k}}}-\frac{2}{t}g_{i{\bar{k}}}\right)Q_{k{\bar{j}}}. \label{Q}
\end{split}
\end{equation}
According to the tensor maximum principle of Hamilton \cite{H1,H2}, to prove $Q_{i{\bar{j}}}\geq 0$, it suffices to show that
\begin{equation}
R_{i{\bar{j}}k{\bar{l}}}Q_{l{\bar{k}}}-\left(\frac{1}{2}+\epsilon\right)\left(R_{p{\bar{j}}}Q_{i{\bar{p}}}+R_{i{\bar{p}}}Q_{p{\bar{j}}}\right) \label{B}
+\left(Q_{i{\bar{k}}}-\frac{2}{t}g_{i{\bar{k}}}\right)Q_{k{\bar{j}}}
\end{equation}
satisfies null-eigenvector condition.
By the well-known result of Bando \cite{B} and Mok \cite{M2}, the nonnegativity of holomorphic bisectional curvature is preserved along the $\epsilon$-K\"{a}hler-Ricci flow. So each term of (\ref{B}) is nonnegative when evaluated at any null-eigenvector of $Q_{i\bar{j}}$. This completes the proof of the theorem.

\begin{center}
\bf Acknowledgements
\end{center}

The first named author would like to thank Professors Q.-K. Lu, S.-K. Wang and K.
Wu for valuable discussions and comments. He would also like to thank Professor H.-W. Xu for helpful suggestions. This work partially supported by the Fundamental Research Funds for Central Universities (Grant No. 2012QNA40).



\end{document}